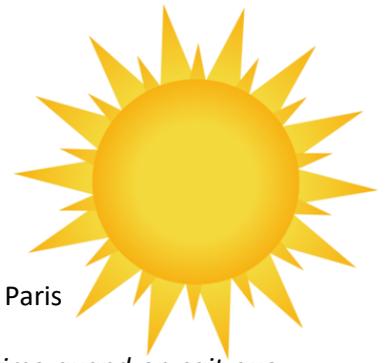

# Le Soleil se lèvera-t-il demain ?

Karim Zayana, Olivier Rioul

Ministère de l'Éducation nationale et de la Jeunesse, Institut Polytechnique de Paris

*Faut-il vivre chaque jour comme s'il était le dernier ? Cette question est légitime quand on sait que Pierre-Simon de Laplace, en 1814, évaluait à 0.99999945 la probabilité que le Soleil se lève à nouveau le lendemain.*

Il n'y a pourtant pas plus trompeur qu'un raisonnement inductif mal construit. D'ailleurs, ne nous met-on pas en garde ? « Ce n'est pas parce que l'on vient de tirer cinq fois 'pile' avec une pièce qu'elle a moins de chances de tomber sur 'pile' au sixième lancer ! » nous avertit le professeur, et d'ajouter en ponctuant ses syllabes « les événements sont in-dé-pen-dants ». À l'inverse, comment ne pas voir une loi des séries quand le destin s'obstine ? Le mathématicien du XX[e] siècle Bertrand Russel aura, lui, cité l'exemple d'une dinde de Noël : généreusement traitée chaque matin, ne prend-elle pas confiance en la main nourricière qui lui tordra pourtant le cou le soir du réveillon ? On peut aussi penser à cette chanson de Brassens :

> On effeuilla vingt fois la marguerite,
> Elle tomba vingt fois sur «pas du tout».
> Et notre pauvre idylle a fait faillite,
> Il est des jours où Cupidon s'en fout.

Qui croire, Laplace, le professeur, Russel ou Brassens ? Et bien chacun d'eux, car tout dépend des hypothèses ! C'est en particulier l'approche bayésienne qui va éclairer la démonstration de Laplace dont nous allons maintenant expliquer les rouages

**Des lancers « indépendants »**

Laplace fonde sa réflexion dès 1774 dans son traité, *Mémoires sur la probabilité des causes par les événements*, et l'illustre 40 ans après sur la succession des jours dans son *Essai philosophique sur les probabilités*. Mais n'allons pas trop vite : voyons d'abord des modèles plus simples.

Commençons par ce cas d'école : la pièce de 1 euro. Supposons-la, dans un premier temps, bien équilibrée, et lançons-la plusieurs fois. Les lancers sont faits dans des conditions qui n'ont rien à voir les unes avec les autres : en langage courant ils sont « indépendants ». On peut donc admettre que les événements de 'pile' ou 'face' successifs seront indépendants au sens mathématique du terme. Comme disait notre professeur : « Ce n'est pas parce que l'on vient de tirer 5 fois 'pile' que cela change les chances du sixième lancer, dont la probabilité de tomber sur 'pile' vaut toujours ½ ». Ainsi la pièce ne garde-t-elle aucune mémoire des lancers précédents.

Ce raisonnement s'applique également à une pièce déséquilibrée, dont la probabilité de tomber sur 'pile' vaudrait $x$ (et donc sur 'face', $1 - x$). Avec des lancers « indépendants », les événements de 'pile' ou 'face' successifs demeurent mathématiquement indépendants : ce n'est pas parce que l'on vient de tirer 5 fois 'pile' que cela change les chances du sixième lancer, dont la probabilité de tomber sur 'pile' vaut toujours $x$.

Compliquons maintenant les choses en prenant deux pièces distinctes, aux probabilités respectives $x_1$ et $x_2$ de tomber sur 'pile' une fois lancées. Pour exagérer, nous pourrions imaginer $x_1$ très voisine de (voire égale à) 1 et $x_2$ très voisine de (voire égale à) 0. Piochons à l'aveugle l'une des deux pièces dans notre porte-monnaie et lançons-la six fois de suite. Dans l'hypothèses où elle tomberait d'abord cinq fois sur 'pile', cela laisse présager qu'on a pioché la première pièce. Cela renforce donc la probabilité qu'elle tombe une sixième fois sur 'pile'. Dans cette nouvelle situation, les lancers restent « indépendants », mais plus les événements de 'pile' ou 'face' successifs ! Mathématiquement, on dit qu'il y a *indépendance conditionnelle* (une fois la pièce choisie), mais *pas indépendance stricte* (inconditionnelle). Notons $E_1, E_2, E_3, E_4, E_5, E_6$ les événements successifs de 'pile' ou 'face' et appelons $N$ le numéro (égal à 1 ou 2) de la pièce choisie au hasard. Le calcul ci-contre prouve que

D'une part, par définition d'une probabilité conditionnelle, avec les notations usuelles,

$$\mathbb{P}_{E_1,\ldots,E_5='\text{pile}'}(E_6 = '\text{pile}') = \frac{\mathbb{P}(E_1,\ldots,E_5,E_6 = '\text{pile}')}{\mathbb{P}(E_1,\ldots,E_5 = '\text{pile}')}$$

D'autre part, par la formule des probabilités totales,

$$\mathbb{P}(E_1,\ldots,E_5 = '\text{pile}')$$
$$= \mathbb{P}_{N=1}(E_1,\ldots,E_5 = '\text{pile}')\mathbb{P}(N = 1)$$
$$+ \mathbb{P}_{N=2}(E_1,\ldots,E_5 = '\text{pile}')\mathbb{P}(N = 2)$$

Convenons de l'équiprobabilité du choix des pièces : $\mathbb{P}(N = 1) = \mathbb{P}(N = 2) = ½$. Puisque les lancers sont indépendants, on obtient le numérateur :

$$\mathbb{P}(E_1,\ldots,E_5 = '\text{pile}') = ½(x_1^5 + x_2^5)$$

Le numérateur s'obtient de la même façon, d'où le résultat.

$$\mathbb{P}_{E_1,\ldots,E_5='\text{pile}'}(E_6 = '\text{pile}') = \frac{x_1^6 + x_2^6}{x_1^5 + x_2^5}$$

Généralisons encore en prenant $k$ pièces distinctes, dont les probabilités de tomber sur 'pile' sont respectivement $x_1, x_2,\ldots, x_k$. Disons que ces pièces ont des probabilités égales d'être choisies dans notre porte-monnaie : $\mathbb{P}(N = 1) = \mathbb{P}(N = 2) = \cdots = \mathbb{P}(N = k) = 1/k$. Tirons une pièce au hasard, puis lançons-la $n+1$ fois. Si les $n$ premiers lancers donnent 'pile',

$$\mathbb{P}_{E_1,\ldots,E_n='\text{pile}'}(E_{n+1} = '\text{pile}') = \frac{\sum_{j=1}^{k} x_j^{n+1}}{\sum_{j=1}^{k} x_j^n} \quad (1)$$

Envisageons enfin des pièces de tailles différentes. Différentes au toucher, leurs probabilités d'être choisies ne sont plus nécessairement égales. Elles valent plus largement $\mathbb{P}(N = 1) = f(x_1)$, $\mathbb{P}(N = 2) = f(x_2),\ldots, \mathbb{P}(N = k) = f(x_k)$ où $f$ désigne une mesure de probabilité. Dès lors,

$$\mathbb{P}_{E_1,\ldots,E_n='\text{pile}'}(E_{n+1} = '\text{pile}') = \frac{\sum_{j=1}^{k} x_j^{n+1} f(x_j)}{\sum_{j=1}^{k} x_j^n f(x_j)} \quad (2)$$

**Tout est en place… pour Laplace !**

Un beau jour de 1814, Laplace se demande avec quelle probabilité le Soleil se lèvera de nouveau le lendemain sachant qu'il s'est toujours levé chaque matin depuis que l'homme est apparu sur Terre. Cette question, devenue célèbre, est reprise dans le programme de mathématiques de spécialité de première en guise d'anecdote historique. Elle peut apparaître incongrue, mais n'oublions pas que le savant, en plus d'être mathématicien, est aussi physicien et philosophe : autant d'angles par lesquels il interroge le monde.

Ainsi, étant posé que le Soleil s'est déjà levé $n$ fois de suite, Laplace se demande s'il reviendra le $n + 1$ ième matin. Il assimile les levers de Soleil à des épreuves « indépendantes » de Bernoulli de même paramètre de succès $p$. Seulement voilà, l'Univers a certes fixé $p$, mais sans le dévoiler. Laplace suppose donc que la valeur de ce paramètre $p$ a été choisi au hasard dans un ensemble de valeurs $x$ comprises entre 0 et 1. Contrairement au cas du porte-monnaie, cet ensemble pourrait comprendre un continuum de valeurs. Ceci amène Laplace à étendre la formule (2) du discret au continu remplaçant les sommes par des intégrales. D'où, en adaptant les notations :

$$\mathbb{P}_{E_1,\ldots,E_n='\text{se lève}'}(E_{n+1}='\text{se lève}') = \frac{\int_0^1 x^{n+1} f(x)\mathrm{d}x}{\int_0^1 x^n f(x)\mathrm{d}x} \tag{3}$$

Un passage à la limite, détaillé en regard, montre que cette probabilité tend vers 1 quand $n$ tend vers l'infini. Cependant, les grands théorèmes de convergence ne seront formalisés que plus tard, au cours du XX[e] siècle par Henri Lebesgue.

Pour son application numérique, Laplace postule *a priori* une densité uniforme (cette hypothèse, devenue classique, est connue sous le nom d' « *a priori* de Laplace »), soit $f = \mathbf{1}_{[0;1]}$. La probabilité recherchée est donc exactement égale à

> Posons $u = x^{n+2}$ dans l'intégrale du numérateur et $u = x^{n+1}$ dans celle du dénominateur. Ces changements de variable conduisent dans un cas à $\mathrm{d}u = (n+2)x^{n+1}\mathrm{d}x$ et dans l'autre à $\mathrm{d}u = (n+1)x^n\mathrm{d}x$. Les bornes d'intégration restent inchangées et le rapport de la formule (3) devient
>
> $$\frac{n+1}{n+2}\frac{\int_0^1 f(u^{1/(n+2)})\mathrm{d}u}{\int_0^1 f(u^{1/(n+1)})\mathrm{d}u}$$
>
> En supposant $f$ continue, les deux intégrales admettent $f(1)$ pour limite par le théorème de convergence dominée. Pourvu que $f(1)$ soit non nul, le rapport étudié équivaut donc au quotient $\frac{n+1}{n+2}$ et tend donc vers 1.

$$\mathbb{P}_{E_1,\ldots,E_n='\text{se lève}'}(E_{n+1}='\text{se lève}') = \frac{\int_0^1 x^{n+1}\,\mathrm{d}x}{\int_0^1 x^n\,\mathrm{d}x} = \frac{n+1}{n+2} \tag{4}$$

Laplace estime ensuite l'âge de l'humanité à 5 000 ans, soit d'après lui, $n = 1\,826\,213$ jours. Ce faisant, il obtient la probabilité annoncée : $1\,826\,214/1\,826\,215 \cong 0.99999945$.

Ce résultat, certes proche de 1 mais pas tout à fait égal à 1, semble contre-intuitif et provient de l'hypothèse *a priori* que l'univers a choisi aléatoirement le paramètre $p$ qui gouvernerait (en partie) le système solaire. L'existence de cet *a priori*, caractéristique de l'approche bayésienne, peut être critiquée. Doit-on pour autant rejeter cette approche ? C'est un débat perpétuel entre bayésiens et fréquentistes !